\documentclass[letterpaper, 10 pt, conference]{ieeeconf}  

\IEEEoverridecommandlockouts                              

\usepackage{dahlin}
\usepackage[latin1]{inputenc}
\usepackage{amsmath}
\usepackage{amssymb}
\usepackage{epsfig}
\usepackage{psfig}
\usepackage{psfrag}
\usepackage{graphics}
\usepackage{color}
\usepackage{pmat}
\usepackage{subfigure}

\newcommand{\Ebf}{{\mathbf E}}

\newcommand{\xbf}{{\mathbf x}}

\newcommand{\zbf}{{\mathbf z}}

\newcommand{\Ccal}{{\mathcal {C}}}

\newcommand{\Ical}{{\mathcal {I}}}

\newcommand{\Pcal}{{\mathcal {P}}}

\newcommand{\Vcal}{{\mathcal {V}}}

\newcommand{\rb}{{\mathbb R}}

\newcommand{\beq}{\begin{equation}}
\newcommand{\eeq}{\end{equation}}

\newcommand{\bnul}{\begin{enumerate}[a)]}
\newcommand{\enul}{\end{enumerate}}

\newcommand{\op}[1]{\operatorname{#1}}

\newcommand{\fin}{\hspace*{\fill}~$\blacksquare$}

\newtheorem{prob}{Problem}
\newtheorem{ejem}{Example}

\title{\LARGE \bf
Particle-based Gaussian process optimization for input design in nonlinear dynamical models
}

\author{Patricio E. Valenzuela, Johan Dahlin, Cristian R. Rojas and Thomas B. Schön
\thanks{This work was supported 
 by the Swedish Research Council under contracts 621-2013-5524 and 621-2009-4017.} 
\thanks{Patricio E. Valenzuela and Cristian R. Rojas are with the Department of Automatic Control and  ACCESS Linnaeus Center, School of Electrical Engineering,
KTH Royal Institute of Technology, SE-100 44 Stockholm, Sweden
(e-mail: {\tt \small \{pva,crro\}@kth.se}).
}%
\thanks{Johan Dahlin is with the Department of Electrical Engineering, Linköping University, SE-581 83 Linköping, Sweden (e-mail: {\tt \small johan.dahlin@liu.se}).}
\thanks{Thomas B. Schön is with the Department of Systems and Control, Uppsala University, SE-751 05 Uppsala, Sweden (e-mail: {\tt \small thomas.schon@it.uu.se}).}
}

\begin{document}

\maketitle
\thispagestyle{empty}
\pagestyle{empty}

\begin{abstract}
We propose a novel approach to input design for identification of nonlinear state space models.
 The optimal input sequence is obtained by maximizing a scalar cost function of the Fisher information matrix. Since the Fisher information matrix is unavailable in closed form, 
  it is estimated using particle methods. In addition, we make use of Gaussian process optimization to find the optimal input and to mitigate the problem of a large computational cost incurred by the particle filter, as the method reduces the number of functional evaluations. 
    Numerical examples are provided to illustrate the performance of the resulting 
     algorithm.
\end{abstract}
\begin{keywords}
System identification, input design, Gaussian process optimization.
\end{keywords}
\section{Introduction}
Input design concerns the maximization of the information retrieved from an experiment. Some of the first contributions in this area have been introduced in \cite{Cox1958,goopay76}. Since then, several approaches to 
 experiment design have been developed (see e.g. \cite{ljung99} 
   and the references therein).


Recently, the problem of input design for the identification of nonlinear dynamical models has gained interest. One of the main difficulties in this case is that a closed form expression for the Fisher information matrix is typically not be available. 
 In addition, the frequency domain techniques employed in the linear case \cite{jansson2005} are no longer valid, which implies that other formulations 
  are required. Contributions in this field consider nonlinear FIR models  \cite{
  larsson2010}, 
 multilevel excitation 
  \cite{forgione2014experiment,decock-gevers-schoukens-13,valenzuela2015graph}, and nonlinear state space models \cite{Valenzuela-Dahlin2014}, among others.

As the Fisher information matrix is unavailable in closed form, we need to 
  rely on estimates. 
   However, such estimates are always subject to uncertainty, which results in difficulties when 
   implementing traditional optimization methods. 

In this work, we explore the reduction of the computational complexity when calculating the objective function used in 
 input design 
  for identification of nonlinear dynamical models. To this end, a Gaussian process optimization (GPO) based algorithm is presented. By assuming that the scalar function of the Fisher information matrix is a realization from a Gaussian process (GP), we can compute its predictive posterior distribution 
   given a set of samples over the feasible 
    set. 
 The predictive posterior distribution 
  acts as a surrogate of the intractable objective function, and is employed to compute the next sample over the feasible set by using an acquisition rule. 
  This technique recursively explores the feasible set to determine the element maximizing a surrogate function. 
   The advantage of this approach when compared with existing techniques is that it can handle uncertainty in the estimates of the objective function, and it drives the exploration of the input space towards those regions where an improvement of the objective function is expected. 

As with most approaches in 
 experiment design, we rely
  on 
  prior information about the system 
   for computing an optimal design. This assumption can be overcome by implementing an adaptive scheme 
 \cite{gerencser2009identification}, 
 or by using a robust input design scheme on top of it \cite{Rojas2007}. However, this 
  is beyond the scope of this paper.


\section{Problem formulation}\label{sec: problem}
Consider the discrete time, nonlinear state space model (SSM) defined for all $t\geq 1$ by
\begin{subequations}
\begin{align}
	x_{t}|x_{t-1}  &\sim  f_\theta (x_{t}|x_{t-1},u_{t-1}), \\
	y_{t}|x_t      &\sim  g_\theta (y_{t}|x_t,u_t), \\
    x_0            &\sim  \mu_\theta(x_0),
\end{align}%
\label{eq:generalSSM}%
\end{subequations}%
where $f_\theta$, $g_\theta$, and $\mu_\theta$ are 
 known 
 probability density functions (pdf) parameterized by the unknown parameter $\theta \in \Theta \subset \rb ^{n_\theta}$. Here, 
  $u_t \in \Ccal \subseteq \rb ^{n_u}$ denotes the input signal, $x_t \in \rb ^{n_x}$ are the (unobserved/latent) internal states, and $y_t \in \rb^{n_y}$ are the measured outputs. In the following, we assume that there exists a $\theta_0 \in \Theta$ such that the pdfs in \eqref{eq:generalSSM} describe the true pdfs of the system when $\theta = \theta_0$, i.e., there is no undermodelling \cite{ljung99}.

The objective 
 is to design 
  $u _{1:T} := (u_{1}, \, \ldots , \, u_{T} )\in \Ccal ^{T}$, such that the parameter $\theta$ in the model \eqref{eq:generalSSM} 
can be identified  with maximum accuracy as defined by 
 a scalar function of the Fisher information matrix $\Ical _F ^{\theta_0}$ \cite{goopay76}, given by 
\beq
\label{eq: prob3}
\Ical _F ^{\theta_0}(u_{1:T}) := \Ebf \left\{ \mathcal{S}(\theta _0) \mathcal{S}^{\top}( \theta_0) | \,  u_{1:T} \right\} \, ,
\eeq
with $\mathcal{S}( \theta_0)$ denoting the score function, i.e., 
\beq
\label{eq: prob2a}
\mathcal{S}(\theta _0) := \left. \nabla \, \ell _\theta (y_{1:T}) \right| _{\theta = \theta _0}\, .
\eeq
Here, $\ell _\theta (y_{1:T})$ 
 denotes the log-likelihood function 
\beq
\label{eq:loglik}
\ell _\theta (y_{1:T}) :=  \log p_\theta (y_{1:T} | u_{1:T}) \,.
\eeq
We note that the expected value in \eqref{eq: prob3} is with respect to the stochastic processes in \eqref{eq:generalSSM}.

In the following, we consider 
 $u_{1:T}$ as 
 a realization of a stationary process. 
  Hence, 
  we will be interested in the 
 \emph{per-sample Fisher information matrix}, given by
\beq
\label{eq: prob3p}
\Ical _F ^{\theta_0, \op{av}}(u_{1:T}) := \dfrac{1}{T} \Ebf _u \left\{ \Ical _F ^{\theta_0}(u_{1:T}) \right\} \, .
\eeq


The input 
 $u _{1:T}$ optimizes a scalar function of \eqref{eq: prob3p}. We define this scalar function as $h \colon \, \rb ^{m \times m} \rightarrow \, \rb$, assumed to be a 
  matrix nondecreasing function \cite[p. 108]{boyvan04}. 

The problem 
 presented here 
 can be summarized as
\begin{prob}\label{prob1}
Find an input signal $u _{1:T} ^{\op{opt}} \in \Ccal ^{T}$ as 
\beq
\label{eq: prob7}
u _{1:T} ^{\op{opt}} := \arg \max _{u _{1:T} \in \Ccal ^{T}} h(\Ical _F^{\theta_0, \op{av}}(u _{1:T})) \, ,
\eeq
where $h \colon \, \rb ^{m \times m} \rightarrow \, \rb$ is a 
 matrix nondecreasing function, and $\Ical _F^{\theta_0, \op{av}}(u _{1:T})$ is given in \eqref{eq: prob3p}.
\fin
\end{prob}
\section{Gaussian process optimization in input design}\label{sec: gpo}
Problem~\ref{prob1} is difficult to solve. One of the main challenges is the characterization of $h(\Ical _F^{\theta_0, \op{av}}(u _{1:T}))$ for all $u_{1:T} \in \Ccal ^T$. Unless assumptions on the model structure \eqref{eq:generalSSM} and the input properties are made, the expression $h(\Ical _F^{\theta_0, \op{av}}(u _{1:T}))$ is often unavailable, and we need to rely on approximations. 
 Moreover, even if an estimate of $h(\Ical _F^{\theta_0, \op{av}}(u _{1:T}))$ is available, part of the existing 
 optimization methods are difficult to implement, since the uncertainty of the estimate is not taken into account. 

Instead, we employ the iterative procedure discussed 
 in \cite{dahlin2014b} to solve Problem~\ref{prob1}. The procedure generates a sequence of iterates $\{ u_{1:T}^{(k)}\}_{k\geq 0}$ for the input excitation. 
 Each iteration consists of three steps: 
\begin{itemize}
\item[(i)] Given 
 $u_{1:T}^{(k)}$, compute an estimate of the objective function $h(\Ical_F ^{\theta_0,\op{av}} (u_{1:T}^{(k)}))$, 
  denoted by $\widehat h_k$. 
\item[(ii)] Given the collection of tuples $\{ u_{1:T}^{(j)}, \, \widehat h_j\}_{j=0}^k$, create a model of the (unavailable) objective function $h(\Ical_F ^{\theta_0,\op{av}} (u_{1:T}))$.
\item[(iii)] Use the model as a surrogate for $h(\Ical_F ^{\theta_0,\op{av}} (u_{1:T}))$ to generate a new iterate $u_{1:T}^{(k+1)}$.
\end{itemize}
The procedure only requires one estimate 
 of $h(\Ical_F^{\theta_0,\op{av}} (u_{1:T}))$ at each iteration, hence keeping the number of estimates 
  as low as possible. 
 Moreover, it requires fewer iterations than a random search, since it focuses on regions of $\Ccal^T$ where an improvement is expected.

For step (i), we employ particle 
 methods to estimate 
  $h (\Ical_F^{\theta_0,\op{av}} (u_{1:T}^{(k)}))$. 
   This 
    is discussed in Section~\ref{sub:if}. 

For steps (ii) and (iii) we use the GPO framework \cite{shahriari2016taking,osborne2009gaussian}. We first
 compute a surrogate of the objective function by modelling it as a Gaussian process, and computing the predictive posterior distribution 
 based on 
 $\{ u_{1:T}^{(j)}, \, \widehat h_j\}_{j=0}^k$.  
  This is discussed in Section~\ref{sub:gpo}.

Then we make use of a heuristic, referred to as the \emph{acquisition rule} (presented in Section~\ref{sub:acq}), to compute $u_{1:T}^{(k+1)}$ based on the GP model. The acquisition rule 
 favours values of $u_{1:T}$ for which the model predicts a large value of the objective function and/or where there is high uncertainty. 
 This 
  establishes a trade-off between exploration and exploitation of the input set. 
    Finally, to employ the GPO framework in input design, we need tractable parameterizations of $\Ccal ^T$, which are discussed 
  in Subsection~\ref{sub:input}. 
\subsection{Estimating the Fisher information matrix}\label{sub:if}
Given $u_{1:T}^{(k)} \in \Ccal ^{T}$, we need to approximate \eqref{eq: prob3p}. 
 To this end, we 
 consider the estimator in 
  \cite{valenzuela2016}, which is based 
   on one estimate of ${\mathcal S}(\theta _0)$ (provided a sufficiently large $T$) to 
    approximate \eqref{eq: prob3p} by \cite{SegalWeinstein1989}
\begin{equation}
\label{eq: MCIfp}
\widehat{ {\Ical}} _F ^{\theta_0,\op{av}} := \dfrac{1}{T}\left[ \sum_{t=1}^T \widehat  {\mathcal S}_t(\theta _0)  (\widehat  {\mathcal S}_t  (\theta _0) )^\top - \dfrac{1}{T} \widehat  {\mathcal S}(\theta _0)  (\widehat  {\mathcal S}  (\theta _0) )^\top \right] \, ,
\end{equation}
where the Fisher identity \cite{CappeMoulinesRyden2005} can be used to write\footnote{For conciseness, we write $\mathbf{v} := v_{1:T}$ for any vector $v_{1:T}$. In addition, we remove the dependence on $k$ of the input, state, and measurements.}
\begin{align}
\label{eq:scoreSSM}
	\mathcal{S}(\theta')
	&=
	\sum_{t=1}^T
	\mathcal{S}_t(\theta')
     \, , \\
\nonumber
\mathcal{S}_t(\theta') &:=
     \dint \nabla \, \xi_{\theta}(x_{t-1:t}) |_{\theta=\theta'}
	p_{\theta'}(x_{t-1:t}| \mathbf{y}, \mathbf{u})
	\dd x_{t-1:t}\, ,
\end{align}
with
\begin{equation*}
	\xi_{\theta}(x_{t-1:t}) :=  \log f_{\theta}( x_{t} |x_{t-1}, u_{t-1} ) + \log
    g_{\theta}( y_{t} |x_t, u_t )\, ,	
\end{equation*}
and 
 $x_{t-1:t}: =\{x_{t-1},x_t\}$.
As we can see from 
 \eqref{eq: MCIfp}, we require an estimate for \eqref{eq:scoreSSM}, 
  which we obtain from 
   particle 
   methods \cite{
LindstenSchon2013}. 

To estimate the score function 
 in \eqref{eq:scoreSSM}, 
 we require the two-step smoothing distribution $p_{\theta}(x_{t-1:t}|\mathbf{y},\mathbf{u})$, which is not available analytically for a general SSM. Instead, we approximate it using an empirical distribution
\begin{align}
	\widehat{p}_{\theta}( \dn x_{t-1:t}|\mathbf{y}, \mathbf{u})
	:=
	\sum_{i=1}^N
	w_t^{(i)}
	\delta_{x_{t-1:t}^{(i)}}
	( \dn x_{t-1:t} ),
	\label{eq:EmpericalTwoStepSmoothingDistribution}
\end{align}
where $x_t^{(i)}$ and $w^{(i)}_t$ denote particle~$i$ and its normalized weight at time $t$. Here, $\{x_t^{(i)},w_t^{(i)}\}_{t=1}^T$ denotes the \textit{particle system} generated by a particle filter and $\delta_{x'}$ denotes the Dirac measure located at $x=x'$. 

\begin{algorithm}[!t]
\caption{\textsf{Bootstrap particle filter (bPF)}}
\small
\textsc{Inputs:} An SSM \eqref{eq:generalSSM}, $\mathbf y$ (observations), $\mathbf u$ (inputs), $N \in {\mathbb N}$ (no. particles).
\\
\textsc{Output:} $\{x_t^{(i)},w_t^{(i)}\}_{i=1}^N$, $t=1,\ldots,T$. 
\algrule[.4pt]
\begin{algorithmic}[1]
	\STATE Sample $x_0^{(i)} \sim \mu_{\theta}(x_0)$ and set $w_0^{(i)}=1/N$.
	\FOR{$t=1$ to $T$}
    \FOR{$i,j=1$ to $N$}
		\STATE
		(\textsf{Resampling})
		Sample $a_t^{(i)}$ from a multinomial distribution with $\mathbf{P}\Big( a_t^{(i)} = j \Big) = w_{t-1}^{(j)}$.
		\STATE
		(\textsf{Propagation})
		Sample $x_t^{(i)} \sim f_{\theta} \Big( x_t^{(i)}\Big| x_{t-1}^{a_t^{(i)}}, u_t \Big)$.
		\STATE Set $x_{0:t}^{(i)} = \Big\{ x_{0:t-1}^{a_t^{(i)}}, x_t^{(i)} \Big\}$.
		\STATE
		(\textsf{Weighting})
		Calculate $\widetilde{w}^{(i)}_t = g_{\theta} \Big( y_t \Big| x_t^{(i)}, u_t \Big)$.
		\STATE Normalize $\widetilde{w}^{(i)}_t$ (over $i$) to obtain $w^{(i)}_t$.
	\ENDFOR
    \ENDFOR
\end{algorithmic}
\label{alg:smc:pf}
\end{algorithm}

\begin{algorithm}[!h]
  \caption{\textsf{Fast forward-filtering backward-simulator with early stopping (fFFBSi-ES)}}
\small
\textsc{Inputs:} Inputs to Algorithm~\ref{alg:smc:pf}, $M \in \mathbb{N}$ (no.\ backward trajectories), $N_{\text{limit}} \in {\mathbb N}$ (limit for when to stop using rejection sampling), $\rho > 0$. \\
\textsc{Output:} $\widehat{\mathcal{I}}_F^{\theta_0,\op{av}}(\mathbf u)$ (estimate of the Fisher information matrix). 
\algrule[.4pt]
\begin{algorithmic}[1]
	\STATE Run Algorithm~\ref{alg:smc:pf} to obtain the particle system $\Big\{x_t^{(i)},w_t^{(i)} \Big\}_{i=1}^N$ for $t=1,\ldots,T$.
    \STATE Sample $\big\{ b_T(j) \big\}_{j=1}^{ M} \sim \textsf{Multi} \big( \{w_T^{(i)} \}_{i=1}^N \big)$.
    \STATE Set $\tilde{x}^{(j)}_T = x_T^{b_T(j)}$ for $j = 1,\ldots,M$.
    \FOR{$t = T-1$ \TO $1$}
    		\STATE $L \gets 1,\ldots,M$.
    		\STATE
	    \COMMENT{Rejection sampling until $N_{\text{limit}}$ trajectories remain.}
    		\WHILE{ $|L| \geq N_{\text{limit}}$ }
    			\STATE $n \leftarrow \textsf{Multi}\big( \{ 1/|L| \}_{i=1}^{|L|} \big)$.
    			\STATE $\delta \leftarrow \emptyset$.
    			\STATE Sample $\big\{ I(k) \big\}_{k = 1}^n \sim \textsf{Multi} \big( \{ w_t^{(i)} \}_{i=1}^N \big)$.
    			\STATE Sample $\big\{ U(k) \big\}_{k = 1}^n \sim \textsf{Uniform}([0,1])$.
    			\FOR{$k = 1$ \TO $n$}
    				\IF{$U(k) \leq f \big( \tilde{x}_{t+1}^{L(k)} | x_{t}^{I(k)} \big) / \rho$}
    					\STATE $b_t( L(k) ) \gets I(k)$.
    					\STATE $\delta \gets \delta \cup \{ L(k) \}$.
    				\ENDIF
    			\ENDFOR
    		\STATE $L \gets L \setminus \delta$.
    \ENDWHILE
    \STATE \COMMENT{Use standard FFBSi for the remaining trajectories \cite{DoucGariverMoulinesOlsson2011}.}
    \FOR {$j \in L$}
    		\STATE Compute $\tilde{w}_{t|T}^{(i,j)} \propto w_t^{(i)} f \big( \tilde{x}_{t+1}^{(j)} | x_{t}^{(i)} \big)$ for $i = 1, \ldots, N$.
    		\STATE Normalize the smoothing weights $\big\{ \tilde{w}_{t|T}^{(i,j)} \big\}_{i=1}^N$.
		\STATE Draw $b_t(j) \sim \mathsf{Multi} \Big(  \big\{ \tilde{w}_{t|T}^{(i,j)} \big\}_{i=1}^N  \Big)$.    		
    \ENDFOR
    \STATE Set $\tilde{x}_{t:T}^{(j)} = \Big\{ x_t^{b_t(j)}, \tilde{x}_{t+1:T}^{(j)} \Big\}$ for $j=1,\ldots,M$.
    \STATE Calculate 
	\vspace{-5mm}
    \begin{align*}
    		\widehat{\mathcal{S}}_t^{(k)}(\theta)
    		=
    		\frac{1}{M} \sum_{j=1}^{M} \nabla \xi_{\theta} \Big( \tilde{x}_{t:t+1}^{(j)} \Big).
    \end{align*}
    \vspace{-5mm}
    \ENDFOR
    \STATE Compute $\widehat{\mathcal{I}}_F^{\theta_0,\op{av}}(u_{1:T})$ using \eqref{eq: MCIfp}.
  \end{algorithmic}
  \label{alg:smc:ffbsi}
\end{algorithm}

Following \cite{valenzuela2016}, here 
 we use the bootstrap particle filter (bPF), see 
   Algorithm~\ref{alg:smc:pf} \cite{DoucetJohansen2011}. However, the estimator 
   \eqref{eq:EmpericalTwoStepSmoothingDistribution} based only on the bPF often suffers from poor accuracy due to particle degeneracy, 
  see e.g.\ \cite{LindstenSchon2013}. 
 To mitigate this problem, we 
 use 
  a particle smoother that introduces a backward sweep after the forward run of the bPF. 
  Here, we use the forward-filtering backwards simulator (FFBSi) with rejection sampling and early stopping 
    \cite{DoucGariverMoulinesOlsson2011}. 

Algorithm~\ref{alg:smc:ffbsi} presents the pseudo-code for the FFBSi. 
 Here, $\mathsf{Multi} (  \{ p^{(i)} \}_{i=1}^N  )$ and $\mathsf{Uniform}([a,b])$ denote the multinomial distribution over $N$ elements, with $p^{(i)}$ being the probability of choosing the $i$-th element, and the uniform distribution with support $[a,b]$, respectively. 
 We note that the parameter $\rho$ required by Algorithm~\ref{alg:smc:ffbsi} is chosen such that $f_{\theta}(x_{t}| x_{t-1}, u_{t-1}) \leq \rho$ for all $t \in \{1,\, \ldots, \, T\}$. The computational complexity of FFBSi is of order $\mathcal{O}(NMT)$, where $N$ and $M$ denote the number of filter and smoother particles, respectively. 
 We refer to \cite{DoucGariverMoulinesOlsson2011} for more details on the effects of $N$, $M$ and $T$ in the accuracy of the 
  estimator. 

\subsection{Modelling the objective function}\label{sub:gpo}
We explore the use of a GP to model the objective function $h(\Ical _F ^{\theta_0,\op{av}}( u_{1:T}))$ 
 \cite{RasmussenWilliams2006}. 
 GPs can be understood as a generalization of the multivariate Gaussian distribution and are commonly used as priors over functions \cite{boyle07}. In this perspective, the posterior obtained by conditioning on 
 the observations corresponds to the functions that could have generated the observations.

In the following, we model the function $h(\Ical_F^{\theta_0,\op{av}} (\cdot))$ as being a priori distributed according to a GP. That is
\begin{equation}
\label{eq:gp}
h(\Ical_F^{\theta_0,\op{av}} (\cdot)) \sim {\mathcal G \mathcal P}\left( m(\cdot) ,\, \kappa (\cdot, \, \cdot) \right) \, ,
\end{equation}
where the process is fully described by the mean function $m(\cdot)$ and the covariance function $\kappa (\cdot , \, \cdot)$. Examples of these functions are a constant for $m$ and a Matérn $s/2$ function for $\kappa$ \cite[p.84]{RasmussenWilliams2006}. 

To simplify the discussion, we will focus on a specific iteration $k$ of the proposed procedure. Let $\mathcal D _k := \{ {\mathbf u}_{1:T} ^{(k)} ,\, \widehat {\mathbf h} _k\}$ 
 denote a set of iterates, where ${\mathbf u}_{1:T} ^{(k)}$ and $\widehat {\mathbf h} _k$ denote matrices obtained by stacking 
  input realizations and estimates of the objective function up to iteration $k$, respectively. In addition, we will assume that
\begin{equation}
\label{eq:dist}
\widehat h_k = h(\Ical _F ^{\theta_0,\op{av}} (u_{1:T}^{(k)})) + z \, ,
\end{equation}
where $z \sim \mathcal N (0, \sigma _z ^2)$, and $\sigma _z > 0$. 
 We note that $\sigma_z$ is unknown a priori, and it needs to be estimated 
  using 
  $\mathcal D _k$. 
The assumption \eqref{eq:dist} seems strict, 
 but the continuous mapping theorem 
 \cite[Theorem~2.7]{billingsley1999}
  shows that 
   the central limit theorem 
    also applies to the estimate $\widehat h_k$, as it is satisfied by 
     \eqref{eq:loglik} asymptotically in the number of particles. 
\begin{ejem}\label{ex:lgss0}
Consider 
\begin{subequations}
\allowdisplaybreaks
\begin{align}
	x_{t+1} | x_t &\sim \mathcal{N} \Big( \phi \, x_t + u_t, 0.1^2 \Big), \\
	y_{t}   | x_t &\sim \mathcal{N} \Big( \alpha\, x_t, 0.1^2 \Big),
\end{align}%
\label{eq:LGSSmodel}%
\end{subequations}%
where the parameters are $\theta=\{\phi,\alpha\}$. We generate $T=10^3$ observations from \eqref{eq:LGSSmodel} 
 with $\theta_0 =\{0.8, 1\}$.

We are interested in estimating $h({\Ical}_F^{\theta_0,\op{av}}(u_{1:T}) ) = \log \det ( {\Ical}_F^{\theta_0,\op{av}}(u_{1:T}) )$, where $u_{1:T}$ is a binary white noise process with values $\{-1,1\}$.
\begin{figure}[t]
  \begin{center}
  \includegraphics[width=0.49\textwidth]{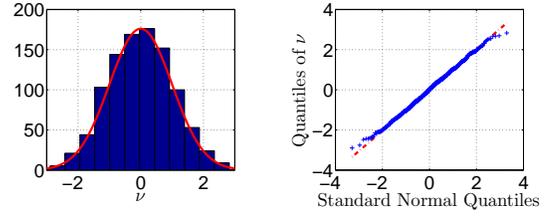}
  \end{center}
  \caption{Left: Histogram of $\nu$ and plot of the scaled pdf of an $\mathcal N (0,1^2)$ distribution (continuous line), Example~\ref{ex:lgss0}. Right: Quantile-quantile plot of the samples of $\nu$ and the $\mathcal N (0,1^2)$ distribution, Example~\ref{ex:lgss0}.}\label{fig:hist}
\end{figure}

The estimate of the Fisher information matrix is obtained using Algorithms~\ref{alg:smc:pf}-\ref{alg:smc:ffbsi}, with $N=2.5 \cdot 10^3$ particles, $M=100$ backward trajectories and $N_{\op{limit}} = \sqrt{N}$ in the fFFBSi smoother. Figure~\ref{fig:hist} shows the histogram based on $10^3$ realizations of the random variable
\beq
\nu := \dfrac{\sqrt{M}(\widehat h - \overline h)}{\sigma_{\sqrt{M} \widehat h}} \, ,
\eeq
where $\widehat h := h({\widehat \Ical}_F^{\theta_0,\op{av}}(u_{1:T}) )$, and $\overline h$, $\sigma _{\sqrt{M} \widehat h}^2$ are 
 the sample mean of $\widehat h$ and variance of $\sqrt{M} \, \widehat h$, respectively. As a comparison, we also present the scaled pdf of an $\mathcal N (0,1^2)$ distribution. We can see that the histogram follows the shape of the pdf of a $\mathcal N (0,1^2)$ distribution. This is also confirmed by the quantile-quantile (QQ) plot in Figure~\ref{fig:hist}, where the quantiles of $\nu$ coincides with those given by an $\mathcal N (0,1^2)$ distribution. \fin 
\end{ejem}

Based on \eqref{eq:dist}, it follows that the predictive posterior distribution is 
\begin{equation}
\label{eq:post}
h(\Ical _F ^{\theta_0,\op{av}}(u_{1:T})) | \mathcal D _k \sim {\mathcal N} \left( \mu (u_{1:T} | \mathcal D _k),\, \sigma^2 (u_{1:T} | \mathcal D _k) + \sigma _z ^2\right) \, ,
\end{equation}
where $\mu (u_{1:T} | \mathcal D _k)$ and $\sigma^2 (u_{1:T} | \mathcal D _k)$ denote the 
 posterior mean and variance given $\mathcal D _k$. From 
  standard results for the Gaussian distribution, we have
\begin{subequations}
\begin{align}
\nonumber
\mu (u_{1:T} | \mathcal D _k) &= m(u_{1:T}) \\
& \hspace{5mm} + \kappa (u_{1:T},\, {\mathbf u}_{1:T}^{(k)}) \Gamma ^{-1} \left\{ \widehat {\mathbf h}_k - m(u_{1:T})\right\}\, , \\
\nonumber
\sigma^2 (u_{1:T} | \mathcal D _k) &=  \kappa (u_{1:T}, \, u_{1:T}) \\
& \hspace{5mm} - \kappa (u_{1:T},\, {\mathbf u}_{1:T}^{(k)}) \Gamma ^{-1} \kappa ({\mathbf u}_{1:T}^{(k)},\, u_{1:T}) \, ,
\end{align}
\label{eq:post2}
\end{subequations}
\hspace{-0.5mm}with $\Gamma :=  \kappa ({\mathbf u}_{1:T}^{(k)} ,\, {\mathbf u}_{1:T}^{(k)}) + \sigma _z ^2 {\mathbf I }_k$, where ${\mathbf I }_k$ denotes the $k\times k$-identity matrix. 

In the GP model introduced here, we use 
 mean and covariance functions that possibly depend on some unknown hyperparameters. In addition, we also need to estimate $\sigma_z$ characterizing the random variable $z$ in \eqref{eq:dist}. To estimate these quantities, we adopt the empirical Bayes procedure, where the marginal likelihood of the data is numerically optimized with respect to the hyperparameters \cite{carlin1996}.

\subsection{Acquisition rules}\label{sub:acq}
To implement step (iii), we need to generate 
 $u_{1:T}^{(k+1)} \in \Ccal ^T$. One option is to perform 
 a random walk over $\Ccal ^T$, which works well provided that the parameterization of $u_{1:T}$ is of small dimension. However, this approach is 
  inefficient as the dimension of the parameterization for $u_{1:T}$ increases.

Instead, we make use of acquisition rules that balance exploration and exploitation of the parameter space and employ 
 the posterior distribution obtained from the GP. Here, we use the expected improvement (EI) technique 
  \cite{jones2001}.

Consider the \emph{predicted improvement}
\begin{equation}
\label{eq:predi}
I(u_{1:T}) := \max \left\{ 0,\, h(\Ical _F ^{\theta_0,\op{av}}(u_{1:T})) -\mu_{\op{max}}-\xi \right\}\, ,
\end{equation}
where $\xi$ is a user defined coefficient balancing exploration and exploitation, and
\begin{equation}
\label{eq:maxu}
\mu_{\op{max}} := \max _{u_{1:T} \in {\mathbf u }_{1:T}^{(k)}} \mu (u_{1:T} | \mathcal D_k)\, ,
\end{equation}
the expected peak of $h(\Ical _F ^{\theta_0,\op{av}} (u_{1:T}))$ at iteration $k$.

By using the posterior distribution obtained from the GP, we define 
 the EI as\footnote{For simplicity, the dependence on $\mathcal D _k$ is dropped from the notation.} 
 \begin{subequations}
\begin{align}
\nonumber
\Ebf \left\{ I(u_{1:T}) \right\} &= \sigma (u_{1:T}) \left\{ Z(u_{1:T}) \Phi (Z(u_{1:T}) ) \right. \\
& \hspace{31mm}\left. - \phi (Z(u_{1:T}))\right\} \, , \\
Z(u_{1:T}) &:= \sigma^{-1} (u_{1:T}) \left\{ \mu (u_{1:T}) -\mu_{\op{max}}-\xi \right\}  ,
\end{align}
\label{eq:eicf}
\end{subequations}
\hspace{-0.9mm}with $\Phi$ and $\phi$ denoting the cumulative distribution function and the pdf 
 of the standard Gaussian distribution, respectively. Then, an acquisition rule is 
\begin{equation}
\label{eq:nextu}
u_{1:T}^{(k+1)} = \arg \max_{u_{1:T} \in \Ccal ^T} \Ebf \left\{ I(u_{1:T}) \left| \mathcal D _k \right. \right\} \, ,
\end{equation}
i.e., the element maximizing the EI. From \eqref{eq:eicf} we see that the EI assigns a large value when both the variance $\sigma (u_{1:T})$ and the mean difference $\mu (u_{1:T}) -\mu _{\op{max}}$ are large, in line with the desired behavior of an acquisition function, as it is explained at the beginning of Section~\ref{sec: gpo}.

\subsection{Parameterizing the input}\label{sub:input}
To implement the GPO for solving the input design problem, we need 
  a 
   parameterization of $\Ccal ^{T}$. 
     Here we briefly explain two options:
\subsubsection{Stationary Markov processes}\label{subsub:graph}
If we restrict 
 $\Ccal$ to be 
  finite 
  and 
   $u_{1:T}$ to be 
    a realization from an $n$-dimensional stationary Markov process of a given order, then the parameterization employed in \cite{valenzuela2015graph} can be used. The parameterization of the input is given by the stationary distribution of the Markov process, which is constrained to 
\begin{multline}
\label{eq: prob6}
\Pcal _\Ccal := \Bigg\{ p_u:\, \Ccal ^{n} \rightarrow \rb \bigg| \, p_u(\xbf) \geq 0 , \, \forall \xbf \in \Ccal ^{n}; \,  
 \\
\left. \sum _ {\xbf \in \Ccal ^{n}} p_u(\xbf) = 1; \right.  \\
 \sum _{v \in \Ccal} p_u(v, \, \zbf) = \sum _{v \in \Ccal} p_u(\zbf , \, v) \, , \forall \zbf \in \Ccal ^{n-1} \Bigg\} \, .
\end{multline}

Following \cite{valenzuela2015graph}, we parameterize \eqref{eq: prob6} as the convex hull of its extreme points, which are computed using graph theoretical techniques. Therefore, the decision variable in this case corresponds to the weighting vector of the extreme points describing an element in $\Pcal _\Ccal$. Assuming that $\Pcal _\Ccal$ has $n_{\Vcal}$ extreme points, then the weighting vector $\alpha := [\begin{matrix} \alpha_1 & \ldots & \alpha _{n_{\Vcal}} \end{matrix}]^\top \in \rb ^{n_{\Vcal}}$ is used to compute 
 $p \in \Pcal _\Ccal$ as
\begin{equation}
\label{eq:convp}
p = \sum_{i=1}^{n_{\Vcal}} \alpha _i p^{(i)} \, ,
\end{equation}
with $\alpha$ satisfying
\begin{subequations}
\begin{align}
\label{eq:alpha1}
\alpha_i &\geq 0 \,, \text{ for all } i \in \{1\, , \ldots ,\, n_{\Vcal}\} \,, \\
\label{eq:alpha2}
\sum_{i=1}^{n_{\Vcal}} \alpha _i &= 1 \, .
\end{align}
\label{eq:constalpha}
\end{subequations}
\noindent In \eqref{eq:convp}, $\{p^{(i)} \}_{i=1}^{n_{\Vcal}}$ corresponds to the probability mass functions (pmf) that are the extreme points of $\Pcal _\Ccal$.

Once a new sample $\alpha \in \rb ^{n_{\Vcal}}$ satisfying\footnote{This can be achieved by sampling $\alpha$ satisfying \eqref{eq:alpha1}, and then normalizing the entries of $\alpha$ to satisfy \eqref{eq:alpha2}.} \eqref{eq:constalpha} is generated, we compute the associated pmf $p \in \Pcal _\Ccal$ by \eqref{eq:convp}, and we generate $u_{1:T}$ by running a Markov chain with stationary distribution $p$.

\subsubsection{Stationary AR processes}
We can restrict $u_{1:T}$ to be a filtered white noise process, as it is proposed in \cite{gopaluni2011}. In this case, the decision variables are the filter coefficients, and the properties of the white noise. For example, we can assume that $u_{1:T}$ is a realization from a stationary AR process 
\begin{equation}
\label{eq:ar}
A(q) \, u_t = e_t \, ,
\end{equation}
where $\{ e_t\}$ is Gaussian white noise, with variance $\sigma_e^2$, and 
\begin{equation}
A(q) := \sum_{i=0}^{n_a} a_i\, q^{-i} \, ,
\end{equation}
with $n_a > 0$ given, $a_i \in \rb$ for all $i \in \{1,\, \ldots,\, n_a\}$, and $a_0 = 1$. For this example, the decision variables are $\sigma_e > 0$, and $\{ a_i \}_{i=1}^{n_a}$, such that $A(q)$ has all its zeros strictly inside the complex unit disc\footnote{This can be guaranteed by factorizing $A(q)$ into first and second order polynomials in $q$, and imposing 
 the constraint on 
  each of these factors.}. 

\subsection{The final procedure}\label{sub:finalproc}
Algorithm~\ref{alg:gpo} presents the resulting procedure for input design using Gaussian process optimization. We note that line~7 introduces a random walk centered at \eqref{eq:nextu} 
 to promote exploration around the expected improvement. 
 We also note that only one functional evaluation is required per iteration, reducing 
 the computational effort when optimizing over $\Ccal^T$.
\begin{algorithm}[!t]
\caption{\textsf{GPO for input design}}
\small
\textsc{Inputs:} Algorithm~\ref{alg:smc:ffbsi}, $K$ (no. iterations) and $u_{1:T}^{(0)} \in \Ccal^T$ (initial excitation).
\\
\textsc{Output:} $\{x_t^{(i)},w_t^{(i)}\}_{i=1}^N$, $t=1,\ldots,T$. 
\algrule[.4pt]
\begin{algorithmic}[1]
	\STATE Sample $u_{1:T}^{(0)}  \in \Ccal ^T$.
	\FOR{$k=0$ to $K$}
		\STATE Use Algorithm~\ref{alg:smc:ffbsi} to compute $\widehat h _k := h(\widehat{\mathcal{I}}_F^{\theta_0,\op{av}}(u_{1:T}^{(k)}))$.
		\STATE Compute \eqref{eq:post}-\eqref{eq:post2} to obtain $h(\Ical _F ^{\theta_0,\op{av}}(u_{1:T})) | \mathcal D _k $.
		\STATE Compute \eqref{eq:maxu} to obtain $\mu_{\op{max}}$.
		\STATE Compute \eqref{eq:nextu} to obtain $\tilde u_{1:T}^{(k+1)}$.
        \STATE Compute $u_{1:T}^{(k+1)}$ as a realization of a random walk centered at $\tilde u_{1:T}^{(k+1)}$.
	\ENDFOR
    \STATE Compute the maximizer of $\mu (u_{1:T} | \mathcal D _K)$ to obtain $u_{1:T}^{\op{opt}}$.
\end{algorithmic}
\label{alg:gpo}
\end{algorithm}

\section{Numerical examples}\label{sec: ex}

\begin{ejem}\label{ex:lgss}
Consider the linear Gaussian state space model in Example~\ref{ex:lgss0}. We are interested in maximizing $h({\Ical}_F^{\theta_0,\op{av}}(u_{1:T}) ) = \log \det ( {\Ical}_F^{\theta_0,\op{av}}(u_{1:T}) )$, where $u_{1:T}$ ($T=10^3$) is a realization 
 of a stationary Markov process (see 
  Section~\ref{sub:input}), with $n_m=1$ and $\Ccal = \{-1,1\}$.

For Algorithm~\ref{alg:gpo}, we use $K = 500$, $\xi = 0.01$, and a random walk centered around the current parametrization of $\tilde u_{1:T}^{(k+1)}$, uniformly distributed on $[-0.01, 0.01]$. The estimate of the Fisher information matrix is obtained using Algorithms~\ref{alg:smc:pf}-\ref{alg:smc:ffbsi}, which are implemented as in Example~\ref{ex:lgss0}. 
 For the prior distribution of $h({\Ical}_F^{\theta_0,\op{av}}(u_{1:T}))$, we consider a constant mean function, and a covariance function composed of a Matérn $s/2$ structure and a constant. The Matérn $s/2$ structure is chosen in this example as it contains information about the smoothness of $h({\Ical}_F^{\theta_0,\op{av}}(u_{1:T}))$. Other choices for the covariance function are also possible and we refer to \cite[Chapter~4]{RasmussenWilliams2006} for more details.

Algorithm~\ref{alg:gpo} is implemented in Matlab using the \verb!fmincon! command for \eqref{eq:nextu} and the \verb+GPML+ toolbox \cite{gpml2015} to infer the hyperparameters and estimate the predictive posterior distribution 
 of $h({\Ical}_F^{\theta_0,\op{av}}(u_{1:T}))$.

The solution obtained from Algorithm~\ref{alg:gpo} is $u_t =1$ for all $t \geq 0$. 
 In this example, a nonzero constant input introduces a nonzero offset in the measurements, which helps to estimate $\theta$ in the presence of process disturbance and measurement noise. As a reference, we draw $u_{1:T}$ as a realization from a binary white noise process with values $\{-1,1\}$. The results are $h({\Ical}_F^{\theta_0,\op{av}}(u_{1:T}^{\op{opt}})) = 14.57$ for the optimal input and $h({\Ical}_F^{\theta_0,\op{av}}(u_{1:T})) = 10.18$ for the binary white noise process. \fin
\end{ejem}
\begin{ejem}\label{ex:1}
Consider the system
\begin{subequations}
\begin{align}
	x_{t+1} | x_t &\sim \mathcal{N} \Big( \frac{1}{\gamma + x^2_t} + u_t, 0.1^2 \Big), \\
	y_{t}   | x_t &\sim \mathcal{N} \Big( \beta x^2_t, 1^2 \Big),
\end{align}%
\label{eq:SGmodel}%
\end{subequations}%
where the parameters are $\theta=\{\gamma,\beta\}$. We generate $T=10^3$ observations from the model with $\theta_0 =\{2,0.8\}$. We note that 
estimating 
 $\gamma$ in \eqref{eq:SGmodel} is inherently difficult, since two different values of $x_t$ can explain $y_t$ equally well.

We consider the same setting and function $h$ 
 as in Example~\ref{ex:lgss}, but we consider 
  three cases for $\Ccal$:
\begin{itemize}
\item Case 1: $\Ccal = \{-1,1\}$.
\item Case 2: $\Ccal = \{-1,0,1\}$.
\item Case 3: $\Ccal = \{-1,-1/3, 1/3,1\}$.
\end{itemize}
\begin{figure}[!b]
  \centering
  \includegraphics[width=0.39\textwidth]{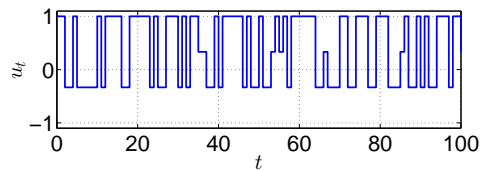}
  \caption{Optimal input $u_{1:T}^{\op{opt}}$ for Case~3 in Example~\ref{ex:1}.}\label{fig:input}
\end{figure}

Table~\ref{tab:results} presents the value of $h^{\op{opt}}:= h(\widehat{\Ical}_F^{\theta_0,\op{av}}(u_{1:T}^{\op{opt}}) )$ for each case, 
 where $u_{1:T}^{\op{opt}}$ corresponds to the optimal input obtained from Algorithm~\ref{alg:gpo}. As comparison, we also compute the value of $h(\widehat{\Ical}_F^{\theta_0,\op{av}}(u_{1:T}) )$, with $\{u_t\}$ binary distributed white noise with values $\{-1,1\}$ (Binary in Table~\ref{tab:results}). 
 We see that the binary white noise process seems to be 
  optimal 
   when $\Ccal = \{-1,1\}$, as it is confirmed by the value of $h^{\op{opt}}$ for Case~1. We also note that adding intermediate values to the input alphabet increases the amount of information in the data, as $h^{\op{opt}}$ is greater in Cases 2 and 3 than in Case~1.
\begin{figure}[t]
  \centering
  \includegraphics[width=0.375\textwidth]{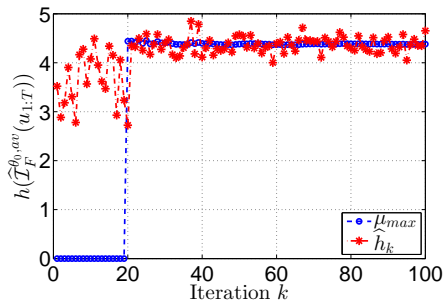}
  \caption{Value of $\widehat h_k$ and $\mu_{\op{max}}$ at iteration $k$ for Case~3 in Example~\ref{ex:1}.}\label{fig:evolution}
\end{figure}

Figure~\ref{fig:input} presents the optimal input obtained for Case~3. We note 
that the optimal input includes a nonzero offset to improve the accuracy of the parameter estimates.

To illustrate the evolution of $\widehat h_k$, we present in Figure~\ref{fig:evolution} the samples $\{ \widehat h_k\}_{k=1}^{100}$, together with the value of $\mu_{\op{max}}$ at every iteration. The first 20 samples are drawn at random from $\Ccal^T$ to provide an initial estimate of the hyperparameters in the GP prior. 
  We note that some of the samples in $\{ \widehat h_k\}_{k=1}^{20}$ are not close to the optimal cost, 
   which is expected due to random sampling. 
  However, once Algorithm~\ref{alg:gpo} is executed from iteration 21 onwards, we observe that the samples are close to $\mu_{\op{max}}$, which implies that the space $\Ccal^T$ is explored only in those regions where $h$ can only increase with respect to the current estimates. Hence, the proposed technique drives the parameter search towards those regions where an improvement in the objective function is expected. \fin
\end{ejem}
\begin{table}[t]
\centering
\caption{$h^{\op{opt}}$ for different input realizations, Example~\ref{ex:1}.}
\begin{tabular}{c||cccc}
  Input & Binary & opt. Case~1 & opt. Case~2 & opt. Case~3 \\ \hline
  $h^{\op{opt}}$ & 4.11 & 4.11 & 4.15 & 4.44 \\
\end{tabular}
\label{tab:results}
\end{table}

\section{Conclusions}\label{sec: conc}
A Gaussian process optimization algorithm for 
 input design for the identification of nonlinear dynamical models has been introduced. The method maximizes a scalar cost function of the Fisher information matrix over the parameter set for the input sequence. Since the objective function is unavailable in closed form, 
  a Gaussian process approach is employed to compute a surrogate function. 
 Numerical examples show that the algorithm can provide a good alternative to solve the input design problem.

Future work on this subject will consider a better estimator of the Fisher information matrix with a better particle smoother, and alternative 
 parameterizations of $\{ u_t\}$. 






\bibliographystyle{IEEEbib}
\bibliography{library,library2a}

\end{document}